\documentclass[12pt]{article}
\usepackage{graphicx}
\usepackage{amssymb, amsmath, amsthm}

\newtheorem{thm}{Theorem}
\newtheorem{lem}[thm]{Lemma}

\begin{document}

\title{Linear complexity of quaternary sequences over $\mathbb{Z}_4$ derived from generalized
cyclotomic  classes modulo $2p$}

\author{Zhixiong Chen$^{1}$, Vladimir Edemskiy$^{2}$\\
1. School of Mathematics, Putian University, \\ Putian, Fujian
351100, P.R. China\\
ptczx@126.com \\
2. Department of Applied Mathematics and Informatics,\\ Novgorod State
University, \\Veliky Novgorod, 173003, Russia \\
Vladimir.Edemsky@novsu.ru}

\maketitle

\begin{abstract}
We determine the exact values of the linear complexity of $2p$-periodic quaternary
sequences over $\mathbb{Z}_4$  (the residue class ring modulo $4$) defined from the generalized cyclotomic
classes modulo $2p$ in terms of the theory of of Galois rings of
characteristic 4, where $p$ is an odd prime. Compared to the case of quaternary
sequences over the finite field of order $4$, it is more difficult and complicated to consider the roots of polynomials in $\mathbb{Z}_4[X]$ due to the zero divisors in $\mathbb{Z}_4$ and hence brings some interesting twists. We prove the main results as follows
$$
\mathrm{Linear ~ Complexity}
=\left\{
\begin{array}{cl}
2p, & \mathrm{if}\,\ p\equiv -3 \pmod 8,\\
2p-1, & \mathrm{if}\,\ p\equiv 3 \pmod 8,\\
p, & \mathrm{if}\,\ p\equiv -1 \pmod {16},\\
p+1, & \mathrm{if}\,\ p\equiv 1 \pmod {16},\\
(p+1)/2, & \mathrm{if}\,\ p\equiv -9 \pmod {16},\\
(p+3)/2, & \mathrm{if}\,\ p\equiv 9 \pmod {16},
\end{array}
\right.
$$
which answers an open problem proposed by Kim, Hong and Song.

\end{abstract}

\section{Introduction}
\label{intro}

Due to applications of quaternary sequences in communication systems, radar and cryptography \cite{GG},
it is of interest to design large families of quaternary sequences.

 Certain quaternary sequences were defined in the literature\footnote{The generalized cyclotomic classes modulo $2p$ are also used to define binary sequences \cite{CDR,DHM,ZZ15}.} by using
generalized cyclotomic classes modulo $2p$ for an odd prime $p$. Let $g$ be an odd number such that $g$  is a primitive root modulo  $p$ and modulo $2p$ simultaneously.
We note that such $g$ always exits, see \cite{IR}. We denote by $\mathbb{Z}_{2p}=\{0,1,\ldots,2p-1\}$ the residue class ring modulo $2p$. Put
$$
D_0=\langle g^{2}\rangle=\{g^{2n} \pmod {2p} : n=0,1,\ldots, (p-3)/2\}\subset \mathbb{Z}_{2p},
$$
and
$$
D_1=\{g^{2n+1} \pmod {2p} : n=0,1,\ldots, (p-3)/2\}\subset \mathbb{Z}_{2p}.
$$
If we write $E_i=\{2u \pmod{2p} : u\in D_i\}$ for $i=0,1$, we have the following partition
$$
  \mathbb{Z}_{2p}= D_0 \cup D_1 \cup E_0 \cup E_1 \cup \{ 0,p \}.
$$
We remark that $D_0 \cup D_1$ is exactly the set of all odd numbers in $\mathbb{Z}_{2p}\setminus \{p\}$
and $E_0 \cup E_1$ is exactly the set of all even numbers in $\mathbb{Z}_{2p}\setminus \{0\}$.

In terms of the generalized cyclotomic classes above,
Chen and Du  \cite{DC} defined a family of
quaternary sequences $(e_u)_{u\ge 0}$ with elements in the finite field  $\mathbb{F}_4=\{0,1, \alpha, 1+\alpha\}$  as
$$
e_u =\left\{
\begin{array}{cl}
0,& \mathrm{if}\,\  i=0 \,\ \mathrm{or}\,\ i \in D_0,\\
1,& \mathrm{if}\,\  i \in D_1,\\
1+\alpha,& \mathrm{if}\,\  i=p \,\ \mathrm{or}\,\ i \in E_0,\\
\alpha,& \mathrm{if}\,\ i \in E_1.
 \end{array}
\right.
$$
They determined the linear complexity of $(e_u)_{u\ge 0}$ in \cite{DC}.
Later Ke, Yang and Zhang \cite{KYZ} calculated their autocorrelation
values\footnote{Ke and Zhang extended to define quaternary cyclotomic sequences of length $2p^m$ \cite{KZ}. Chang and Li defined quaternary cyclotomic sequences of length $2pq$ \cite{CL}. Both are over the finite field $\mathbb{F}_4$.}
. In fact, before \cite{DC,KYZ} Kim, Hong and Song \cite{KHS}
defined another family of quaternary sequences $(s_u)_{u\ge 0}$ with
elements in $\mathbb{Z}_4=\{0,1, 2, 3\}$, the residue class ring
modulo $4$, as follows
 \begin{equation}
\label{eq1}
s_u=\left\{
\begin{array}{cl}
0,& \mathrm{if}\,\  i=0 \,\ \mathrm{or}\,\ i \in D_0,\\
1,& \mathrm{if}\,\  i \in D_1,\\
2,& \mathrm{if}\,\  i=p \,\ \mathrm{or}\,\ i \in E_0,\\
3,& \mathrm{if}\,\ i \in E_1.
 \end{array}
\right.
\end{equation}
They derived the periodic\footnote{The period of $(s_u)_{u\ge 0}$ is $2p$.} autocorrelation function of $(s_u)_{u\ge 0}$. However, as we know, the linear complexity of $(s_u)_{u\ge 0}$ is still open  since
 it faces more difficulties due to the  phenomenon of zero divisors in $\mathbb{Z}_4$. In this work, we will develop a way to solve this problem using the theory of Galois rings of
characteristic 4. We note that there are many quaternary sequences
over $\mathbb{Z}_4$ have been investigated in the literature, see
e.g., \cite{US96,US98,EI1,EI2,EI3}.

We recall that the
\emph{linear complexity} $LC((s_u)_{u\ge 0})$  of $(s_u)_{u\ge 0}$ above  is the least order $L$ of a linear
recurrence relation (i.e., linear feedback shift register, or LFSR for short) over $\mathbb{Z}_4$
$$
s_{u+L} + c_{1}s_{u+L-1} +\ldots +c_{L-1}s_{u+1}+ c_Ls_u=0\quad
\mathrm{for}\,\ u \geq 0,
$$
which is satisfied by $(s_u)_{u\ge 0}$ and where $c_1, c_2, \ldots,
c_{L}\in  \mathbb{Z}_4$, see \cite{US00}. The \emph{connection polynomial} is $C(X)$ given by $1+c_1X+\ldots+c_LX^L$.
We note that $C(0)=1$. Let
$$S(X)=s_0+s_1X+\ldots+s_{2p-1}X^{2p-1}\in  \mathbb{Z}_4[X]$$
be the \emph{generating polynomial} of $(s_u)_{u\ge 0}$. Then an
LFSR with a connection polynomial $C(X)$ generates $(s_u)_{u\ge 0}$,
if and only if \cite{US00},
$$
S(X)C(X) \equiv 0 \pmod {X^{2p}-1}.
$$
That is,
 \begin{equation}
\label{eq2}
LC((s_u)_{u\ge 0})=\min\{\deg(C(X)) : S(X)C(X) \equiv 0 \pmod {X^{2p}-1}\}.
\end{equation}

Let $r$ be the order of 2 modulo $p$. We denote by  $GR(4^{r},4)$ the Galois ring of order $4^r$ of characteristic 4,
which is isomorphic to the residue  class ring $\mathbb{Z}_4[X]/(f(X))$, where $f(X)\in \mathbb{Z}_4[X]$
is a \emph{basic irreducible polynomial} of degree $r$, see \cite{W1,M}.
 The group of units of $GR(4^r,4)$ is denoted by $GR^*(4^r,4)$, which contains a cyclic subgroup of order $2^r-1$.
 Since $p|(2^r-1)$, let $\beta\in GR^*(4^r,4)$ be of order $p$. Then we find that\footnote{In the context we always suppose that $\gamma\in GR^*(4^r,4)$ is of order $2p$.} $\gamma=3\beta\in GR^*(4^r,4)$ is of order $2p$.
From (\ref{eq2}), we will consider the values $S(\gamma^v)$ for
$v=0,1, \ldots, 2p-1$, which allow us to derive the linear
complexity of $(s_u)_{u\ge 0}$. Due to  $S(X)\in \mathbb{Z}_4[X]$,
we cannot consider it in the same way as those in finite fields. For
example, 1  and  3 are the roots of $2X-2\in \mathbb{Z}_4[X]$, but
$2X-2$ is not divisible by $(X-1)(X-3)$, i.e., in the ring
$\mathbb{Z}_4[X]$ the number of roots of a polynomial can be greater
than its degree. So we need to develop some necessary technique
here. Indeed, the theory of Galois ring enters into our problem by
means of the following lemmas.

\begin{lem}
\label{pre-root} Let $P(X)\in \mathbb{Z}_4[X]$ be a non-constant
polynomial.  If $\xi\in GR(4^r,4)$ is a root of $P(X)$, we have
$P(X)=(X-\xi)Q_1(X)$ for some polynomial $Q_1(X)\in GR(4^r,4)[X]$.

 Furthermore, if $\eta\in GR(4^r,4)$ is another root of $P(X)$ and $\eta-\xi$ is a unit, we have $P(X)=(X-\xi)(X-\eta)Q_2(X)$, where $Q_1(X)=(X-\eta)Q_2(X)$.
\end{lem}

\begin{lem}
\label{root}
  Let $\gamma\in GR^*(4^r,4)$ be of order $2p$, and let
$P(X)\in \mathbb{Z}_4[X]$ be any non-constant polynomial.

(I). If $P(\gamma^v)=0$ for all $v \in D_i$, where $i=0,1$, then we
have
$$
P(X)=P_1(X) \prod_{v \in D_i} (X-\gamma^v)
$$
for some polynomial $P_1(X) \in  GR(4^r,4) [X]$. Similarly, If
$P(\gamma^v)=0$ for all $v \in E_i$, where $i=0,1$, then we have
$$
P(X)=P_2(X) \prod_{v \in E_i} (X-\gamma^v)
$$
for some polynomial $ P_2(X) \in  GR(4^r,4)[X]$.

(II). If $P(\gamma^v)=0$ for all $v \in \{p\} \cup D_0 \cup D_1$, then we have
$$
P(X)=P_3(X)(X^p+1)
$$
for some polynomial $ P_3(X) \in GR(4^r,4)[X]$.
 Similarly, if $P(\gamma^v)=0$ for all $v \in \{0\}
\cup E_0 \cup E_1$, then we have
$$
P(X)=P_4(X)(X^p-1)
$$
for some polynomial $ P_4(X) \in \mathbb{Z}_4[X]$.

(III). If $P(0)=1$, $P(\gamma^v)=0$ for $v\in \mathbb{Z}_{2p}\setminus \{0,p\}$, and $P(\pm 1)\in \{0,2\}$,  then we have $\deg P(X)\geq 2p-1$. Furthermore, if either $P(1)=P(-1)=0$ or $P(1)=P(-1)=2$, we have $\deg P(X)\geq 2p$.
\end{lem}

We give a proof of both lemmas in the Appendix for the convenience of the reader.

\section{Linear complexity of $(s_u)_{u\ge 0}$}
\label{sec:1}

\subsection{Auxiliary lemmas}

We describe a relationship among $D_0,D_1,E_0$ and $E_1$.

\begin{lem}\label{relationD}
Let $i,j\in \{0,1\}$.

(I). For $v\in D_i$, we have
$$vD_j\triangleq \{vu \pmod{2p} : u\in D_j\}=D_{i+j \bmod 2},$$
and
$$v E_j\triangleq \{vu \pmod{2p} : u\in E_j\}=E_{i+j \bmod 2}.$$

(II). For $v\in E_i$, we have
$$vD_j\triangleq \{vu \pmod{2p} : u\in D_j\}=E_{i+j \bmod 2},$$
and
$$v E_j\triangleq \{vu \pmod{2p} : u\in E_j\}=E_{i+j \bmod 2}$$
if $p\equiv \pm 1\pmod {8}$, and otherwise
$$v E_j\triangleq \{vu \pmod{2p} : u\in E_j\}=E_{i+j+1 \bmod 2}.$$

(III). If $p\equiv \pm 1\pmod {8}$, we have
$$
D_i=\{(v+p)\pmod {2p} : v\in E_i\},
$$
and otherwise
$$
D_{i+1\bmod 2}=\{(v+p)\pmod {2p} : v\in E_i\}.
$$

(IV). If $p\equiv \pm 1\pmod {8}$, we have
$$E_i=\{(v+p)\pmod {2p} : v\in D_i\},$$
and otherwise
$$E_{i+1\bmod 2}=\{(v+p)\pmod {2p} : v\in D_i\}.$$

(V). If $p\equiv \pm 1\pmod {8}$, we have
$$
E_i=\{u+p : u\in D_i, u<p\} \cup \{u-p : u\in D_i, u>p\}.
$$
\end{lem}
Proof. (I). If $v\in D_i$ for $i=0,1$ and $u\in D_j$ for $j=0,1$
then we can write that $v\equiv g^{i+2k} \pmod {2p}, 0\leq k \leq
(p-3)/2$ and $u\equiv g^{j+2l} \pmod {2p}, 0\leq l \leq (p-3)/2$.
So, $vu\equiv g^{i+j+2k+2l} \pmod {2p}$, hence $vu \in D_{i+j \bmod
2}.$ Since $|vD_i|=|D_{i+j \bmod 2}|$, it follows that  $vD_i=D_{i+j
\bmod 2}$. The equality $v E_j=E_{i+j+1 \bmod 2}$ may be proved
similarly as the first.

(II). Let $v\in E_i$. We write  $v\equiv 2u \pmod {2p}, u\in D_i$.
Therefore, by (I) and  our definitions we have that
$$
vD_j=2uD_j=2D_{i+j \bmod 2}=E_{i+j \bmod 2}.
$$

Now, we consider $vE_j$.  First, we have by (I) again
$$v E_j=2u E_j=2E_{i+j \bmod 2}.$$
Second, for any $w\in 2E_{i+j \bmod 2}$, we can write $w\equiv 4a \pmod {2p}$ for $a\in D_{i+j \bmod 2}$.
Clearly $w$ is even and $w\in E_0\cup E_1$, so we have $w\equiv 2b \pmod {2p}$ for $b\in D_0\cup D_1$. Then we have
$b\equiv 2a \pmod {p}$.

For $p\equiv \pm 1 \pmod 8$, in which case $2$ is a quadratic
residue modulo $p$ \cite{IR}, we have $b\in D_{i+j \bmod 2}$, which
leads to $w\equiv 2b \pmod {2p}\in E_{i+j \bmod 2}$, i.e., $2E_{i+j
\bmod 2}\subseteq E_{i+j \bmod 2}$. Since $2E_{i+j \bmod 2}$ and
$E_{i+j \bmod 2}$ have the same cardinality, it follows that $v
E_j=2E_{i+j \bmod 2}= E_{i+j \bmod 2}$.

The case of $p\equiv \pm 3 \pmod 8$ follows in a similar way, in
which case $2$ is a quadratic non-residue modulo $p$.

(III).  Let $p\equiv \pm 1 \pmod 8$, in which case $2$ is a
quadratic residue modulo $p$  \cite{IR}. Then we can find when
$\overline{v}$ runs through $D_0$ (resp. $D_1$), $p+2\overline{v}$
modulo $2p$ runs through $D_0$ (resp. $D_1$). Since otherwise, if
$p+2\overline{v}_0 \pmod {2p} \in D_1$ for some $\overline{v}_0\in
D_1$, then  we write $p+2\overline{v}_0\equiv g^{1+2k_0} \pmod {2p}$
for some integer $k_0$, from which we derive $2\overline{v}_0\equiv
g^{1+2k_0} \pmod {p}$. It leads to the result that  $2$ is a
quadratic non-residue modulo $p$, a contradiction. So,
$D_i=\{(v+p)\pmod {2p} : v\in E_i\} $ if $p\equiv \pm 1\pmod {8}$.

The equality $D_{i+1\bmod 2}=\{(v+p)\pmod {2p} : v\in E_i\}$ for
$p\equiv \pm 3\pmod {8}$ is proved similarly as the first. Here,
if $\overline{v}$ runs through $D_0$ (resp. $D_1$),then
$p+2\overline{v}$ modulo $2p$ runs through $D_1$ (resp. $D_0$).

(IV) Comes from (III).

(V).
In fact first, the set $\{u+p : u\in D_1, u<p\} \cup \{u-p : u\in D_1, u>p\}$ exactly contains $|D_1|$ many even numbers.
Second, we suppose that $a\in E_0$ for some $a\in \{u+p : u\in D_1, u<p\} \cup \{u-p : u\in D_1, u>p\}$.
 Write $a\equiv 2v \pmod {2p}$ for some  $v\in D_0$. From the definition of $D_0$, we see that $v$ is a quadratic residue modulo $p$. Then $a$ is a quadratic residue modulo $p$ due to $p\equiv \pm 1\pmod {8}$, in which case
$2$ is a quadratic residue modulo $p$ \cite{IR}. However, $a$ is of the form $u+p$ or $u-p$ for some $u\in D_1$, and
$a\equiv u\pmod p$ is a quadratic non-residue modulo $p$, a contradiction. So $\{u+p : u\in D_1, u<p\} \cup \{u-p : u\in D_1, u>p\}\subseteq E_1$, and both have the same cardinality. \qed\\

For $i=0,1$, let
$$
S_i(X)=\sum_{u \in D_i} X^u
$$
  and
$$
T_{i}(X)=S_i(X^2)=\sum_{u \in E_i} X^u \pmod{X^{2p}-1}.
$$
According to (\ref{eq1}), the generating polynomial of $(s_u)_{u\ge
0}$ is
 \begin{equation}
\label{SX}
S(X)=2X^p+S_1(X)+2T_0(X)+3T_1(X).
\end{equation}
As mentioned before, we will consider the values $S(\gamma^v)$ for a
unit $\gamma\in GR^*(4^r,4)$ of order $2p$ and $v=0,1, \ldots,
2p-1$. According to the definitions of $D_0,D_1,E_0$ and $E_1$, we
will describe $S(\gamma^v)$ in the following Lemma in terms of
$S_0(\gamma)$ (or $S_1(\gamma)$) due to the fact that in the ring
$GR(4^{r},4)$
\begin{equation}
\label{eq3}
S_0(\gamma)+S_1(\gamma)=\sum\limits_{u\in D_0\cup
D_1}\gamma^u=\sum\limits_{j=0}^{p-1}\gamma^{2j+1}-\gamma^p=0+1=1.
\end{equation}

\begin{lem}
\label{l1}
Let $\gamma\in GR^*(4^r,4)$ be of order $2p$, and let S(X) be the generating polynomial of $(s_u)_{u\ge
0}$ described in (\ref{SX}).

(I). If $p\equiv \pm 3 \pmod 8$, we have
$$
 S(\gamma^v)=\left\{
\begin{array}{cl}
1-2S_0(\gamma),& \mathrm{if}\,\  v \in D_0,\\
-1+2S_0(\gamma),& \mathrm{if}\,\  v \in D_1,\\
3,&\mathrm{if}\,\ v \in E_0\cup E_1.
 \end{array}
\right.
$$

(II). If $p\equiv \pm 1 \pmod 8$, we have
$$
 S(\gamma^v)=\left\{
\begin{array}{cl}
0,& \mathrm{if}\,\ v \in D_0\cup D_1,\\
2-2S_0(\gamma),& \mathrm{if}\,\  v \in E_0,\\
2S_0(\gamma),&\mathrm{if}\,\ v \in E_1.
 \end{array}
\right.
$$
\end{lem}
Proof. (I). Let $p\equiv \pm 3 \pmod 8$. By Lemma \ref{relationD}(I) we first get
$$
 S_1(\gamma^v)=\left\{
\begin{array}{cl}
1-S_0(\gamma),& \mathrm{if}\,\ v \in D_0,\\
S_0(\gamma),& \mathrm{if}\,\ v \in  D_1.
\end{array}
\right.
$$
Second, for any $v\in E_j$ for $j\in\{0,1\}$, write $v=2\overline{v}$ for $\overline{v}\in D_j$. We see that
$p+2\overline{v}\in D_{j+1}$ by Lemma \ref{relationD}(III) and $\gamma^{v}=\gamma^{2\overline{v}}=-\gamma^{p+2\overline{v}}$, by Lemma \ref{relationD}(I) we derive
$$
 S_1(\gamma^v)= S_1(-\gamma^{p+2\overline{v}})=-\sum\limits_{u\in D_1}\gamma^{u(p+2\overline{v})}
 =\left\{
\begin{array}{cl}
-\sum\limits_{w\in D_0}\gamma^{w},& \mathrm{if}\,\ \overline{v} \in D_0,\\
-\sum\limits_{w\in D_1}\gamma^{w},& \mathrm{if}\,\ \overline{v} \in  D_1,
 \end{array}
\right.
$$
which leads to
 $$
 S_1(\gamma^v)=\left\{
\begin{array}{cl}
-S_0(\gamma),& \mathrm{if}\,\  v \in E_0,\\
-1+S_0(\gamma),&\mathrm{if}\,\ v \in E_1.
 \end{array}
\right.
$$
Similarly, by Lemma \ref{relationD}(I)-(IV), we have
$$
 T_0(\gamma^v)=\left\{
\begin{array}{cl}
-1+S_0(\gamma),& \mathrm{if}\,\ v \in D_0,\\
-S_0(\gamma),& \mathrm{if}\,\ v \in  D_1,\\
S_0(\gamma),& \mathrm{if}\,\  v \in E_0,\\
1-S_0(\gamma),&\mathrm{if}\,\ v \in E_1,
 \end{array}
\right.
$$
and
$$
 T_1(\gamma^v)=\left\{
\begin{array}{cl}
-S_0(\gamma),& \mathrm{if}\,\ v \in D_0,\\
-1+S_0(\gamma),& \mathrm{if}\,\ v \in  D_1,\\
1-S_0(\gamma),& \mathrm{if}\,\  v \in E_0,\\
S_0(\gamma),&\mathrm{if}\,\ v \in E_1.
 \end{array}
\right.
$$
Then putting everything together, we get the first assertion.

 The second assertion of this lemma can be  proved in a similar way. \qed

So in order to determine the values of $S(\gamma^v)$, it is
sufficient to calculate $S_0(\gamma)$. We need the parameter $[i,j]$ for $i,j\in \{0,1\}$, which is the cardinality of the set $(1+D_i)\cap E_j$, i.e.,
$$
[i,j]=|(1+D_i)\cap E_j|,
$$
where $1+D_i=\{1+u \pmod {2p} : u\in D_i\}$.
\begin{lem}\label{lh}
With notations as before. We have
$$
[0,0]=
\left\{
\begin{array}{cl}
(p-5)/4,& \mathrm{if}\,\ p\equiv  1 \pmod {8},\\
(p-3)/4,& \mathrm{if}\,\ p\equiv  7 \pmod {8},\\
(p-1)/4,& \mathrm{if}\,\  p\equiv 5 \pmod {8},\\
(p+1)/4,& \mathrm{if}\,\  p\equiv 3 \pmod {8},
 \end{array}
\right.
$$
and
$$
[0,1]=
\left\{
\begin{array}{cl}
(p-1)/4,& \mathrm{if}\,\ p\equiv  1 \pmod {8},\\
(p+1)/4,& \mathrm{if}\,\ p\equiv  7 \pmod {8},\\
(p-5)/4,& \mathrm{if}\,\  p\equiv 5 \pmod {8},\\
(p-3)/4,& \mathrm{if}\,\  p\equiv 3 \pmod {8},\\
 \end{array}
\right.
$$
\end{lem}
Proof. Since $g$ used above is also a primitive modulo $p$, we write
$$
H_0=\{g^{2n} \pmod p : n=0,1,\ldots,(p-3)/2\}
$$
and
$$
H_1=\{g^{1+2n} \pmod p : n=0,1,\ldots,(p-3)/2\}.
$$
We find that for $i=0,1$
$$
\{u \pmod p : u\in D_i\}=H_i
$$
and
$$
\{2u \pmod p : u\in D_i\}=H_{i+\ell\bmod 2},
$$
where $\ell = 0$ if $p\equiv \pm 1 \pmod 8$ and $\ell =1$ if $p\equiv \pm 3 \pmod 8$,
  i.e., $\ell=0$ if 2 is a quadratic residue modulo $p$, and $\ell=1$ otherwise \cite{IR}.
Therefore,
$$
[i,j]=|(1+D_i)\cap E_j|=|(1+H_i)\cap H_{j+\ell\bmod 2}|.
$$
We conclude the proof by applying the values of $|(1+H_i)\cap H_{j}|$ computed in \cite{H}.
 \qed

With the values of $[0,0]$ and $[0,1]$, we prove the following statement, which is a generalization of
\cite[Theorem 1]{E}.

\begin{lem}
\label{l3}
 Let $\gamma\in GR^*(4^r,4)$ be of order $2p$.
Then we have
$$
(S_0(\gamma))^2=S_0(\gamma)+ \left\{
\begin{array}{cl}
(p-1)/4,& \mathrm{if}\,\ p\equiv  1 \pmod {8},\\
(p+1)/4,& \mathrm{if}\,\  p\equiv -1 \pmod {8},\\
(p+1)/4,& \mathrm{if}\,\ p\equiv  3 \pmod {8},\\
(p-1)/4,& \mathrm{if}\,\  p\equiv -3 \pmod {8}.
 \end{array}
\right.
$$
\end{lem}
Proof.  By the definition of $S_0(X)$  we have
$$
(S_0(\gamma))^2=\sum_{l,m=0}^{(p-3)/2} \gamma^{g^{2l}+g^{2m}}=\sum_{l,m=0}^{(p-3)/2}
\gamma^{g^{2m}(g^{2(l-m)}+1)}.
$$
For each fixed $m$, since the order of $g$ modulo $2p$ is $p-1$, we see that $l-m$ modulo $(p-1)$ runs through the range $0,1,\ldots,(p-3)/2$ if $l$ does. So we have
\begin{equation}\label{eq4}
(S_0(\gamma))^2=\sum_{m,n=0}^{(p-3)/2}
\gamma^{g^{2m}(g^{2n}+1)}.
\end{equation}

Since $g$ is odd, we see that $g^{2n}+1 \pmod {2p}$ is even for any $n$. That is,
$g^{2n}+1 \pmod {2p} \in E_0\cup E_1\cup \{0\}$.
So we consider
$g^{2n}+1 \pmod {2p}$ in three different cases.

Case 1. Let
$$
N_0=\{n : 0\le n\le (p-3)/2,~ g^{2n}+1 \pmod {2p} \in E_0\}.
$$
In fact, the cardinality $|N_0|$ of $N_0$ equals $[0,0]$.
For each $n \in N_0$, as the proof of Lemma \ref{l1}
we obtain that by (\ref{eq3})
\begin{eqnarray*}
\sum_{m=0}^{(p-3)/2}
\gamma^{g^{2m}(g^{2n}+1)}&=&\sum_{v\in D_0}
\gamma^{2v}=S_0(\gamma^{2})=S_0(-\gamma^{p+2})\\
&=&\left\{
\begin{array}{cl}
-S_0(\gamma),& \mathrm{if}\,\ p\equiv  \pm 1 \pmod {8},\\
-1+S_0(\gamma),& \mathrm{if}\,\  p\equiv  \pm 3 \pmod {8}.
 \end{array}
\right.
\end{eqnarray*}

Case 2. Similar to Case 1, we let
$$
N_1=\{n : 0\le n\le (p-3)/2,~ g^{2n}+1 \pmod {2p} \in E_1\}.
$$
Then the cardinality $|N_1|$ equals $[0,1]$.
Now for each $n \in N_1$,
we obtain that
\begin{eqnarray*}
\sum_{m=0}^{(p-3)/2}
\gamma^{g^{2m}(g^{2n}+1)}&=&\sum_{v\in D_1}
\gamma^{2v}=S_1(\gamma^{2})=S_1(-\gamma^{p+2})\\
&=&\left\{
\begin{array}{cl}
-1+S_0(\gamma),& \mathrm{if}\,\ p\equiv  \pm 1 \pmod {8},\\
-S_0(\gamma),& \mathrm{if}\,\  p\equiv  \pm 3 \pmod {8}.
 \end{array}
\right.
\end{eqnarray*}

Case 3. There is an $n$ such that $(g^{2n}+1)\equiv 0\pmod {2p}$ if and only if
$p\equiv  1 \pmod {4}$. In this case, we have $n=(p-1)/4$ and $\sum_{m=0}^{(p-3)/2}
\gamma^{g^{2m}(g^{2n}+1)}=(p-1)/2$.

Let  $p\equiv 1 \pmod 8$. Using \eqref{eq4} we obtain that
$$
(S_0(\gamma))^2=|N_0|\cdot (-S_0(\gamma))+|N_1|\cdot (-1+S_0(\gamma))+
(p-1)/2.
$$
Then we get the desired result by using the values of$[0,0]$ (=$|N_0|$) and $[0,1]$ ($=|N_1|$) in Lemma \ref{lh}.

The  assertions for
$p\equiv -1, 3, -3 \pmod 8$ can be obtained in a similar way.
 \qed

With the help of Lemma \ref{l3} we now deduce
the values of
  $S_0(\gamma)$. It is clear that $S_0(\gamma) \in GR^*(4^r,4)$ or
$S_1(\gamma) \in GR^*(4^r,4)$ from (\ref{eq3}). Therefore, without loss of
generality we always suppose that $S_0(\gamma) \in GR^*(4^r,4)$. (Of course, if one supposes that  $S_1(\gamma) \in GR^*(4^r,4)$, then  $S_1(\gamma)$ will be used in the context.)

\begin{lem}
\label{l4}
 Let $\gamma\in GR^*(4^r,4)$ be of order $2p$ with  $S_0(\gamma)=\sum_{u \in D_0} \gamma^u \in GR^*(4^r,4)$. We have
$$
S_0(\gamma)=
\left\{
\begin{array}{cl}
1,& \mathrm{if}\,\ p\equiv \pm 1 \pmod {16},\\
\rho,& \mathrm{if}\,\  p\equiv \pm 5 \pmod {16},\\
3,&\mathrm{if}\,\ p\equiv \pm 9 \pmod {16},\\
2+\rho,& \mathrm{if}\,\  p\equiv \pm 13 \pmod {16},
 \end{array}
\right.
$$
where $\rho$ satisfies the equation
$\rho^2+3\rho+3=0$ over $\mathbb{Z}_4$.
\end{lem}
\emph{Proof} Let $p\equiv \pm 1 \pmod {16}$.   Then, by Lemma
\ref{l3}, we obtain that $ (S_0(\gamma))^2=S_0(\gamma)$. Under given
assumptions about $S_0(\gamma)$, we have $S_0(\gamma)=1.$ The other
assertions of this lemma can be  proved in a similar way.
 \qed

\begin{lem}
\label{l5} Let $\gamma\in GR^*(4^r,4)$ be of order $2p$ with  $S_0(\gamma)=\sum_{u \in D_0} \gamma^u \in GR^*(4^r,4)$, and let S(X) be the generating polynomial of $(s_u)_{u\ge
0}$ described in (\ref{SX}).

\noindent
(I). For any odd prime $p$, we have
$$
S(\gamma^v)=\left\{
\begin{array}{cl}
p+1,&\mathrm{if}\,\   v=0,\\
2,&\mathrm{if}\,\   v=p.
\end{array}
\right.
$$
(II). If $p\equiv \pm 3 \pmod 8$, we have
$$
S(\gamma^v)\in GR^*(4^r,4), ~ \mathrm{for~ all} ~ v\in D_0\cup D_1\cup E_0\cup E_1,
$$
(III). If $p\equiv \pm 1 \pmod 8$, we have
$$
S(\gamma^v)=\left\{
\begin{array}{cl}
0,& \mathrm{if}\,\  v\in D_0\cup D_1\cup E_0,\\
2,&\mathrm{if}\,\   v\in E_1.
\end{array}
\right.
$$
\end{lem}
Proof. (I) can be checked easily. (II) and (III) follow  immediately from Lemmas \ref{l1} and \ref{l4}. \qed

In the following subsections, we will derive linear complexity of $(s_u)_{u\ge 0}$ in
\eqref{eq2} by considering the
factorization of $S(X)$.

\subsection{Linear complexity for the case $p\equiv \pm
3 \pmod 8$}
\label{subsec:3}

 \begin{thm}
 \label{t1}
  Let $(s_u)_{u\ge 0}$ be the quaternary sequence over $\mathbb{Z}_4$ defined by \eqref{eq1}. Then the  linear complexity of $(s_u)_{u\ge 0}$
  satisfies
$$
LC((s_u)_{u\ge 0})
=\left\{
\begin{array}{cl}
2p, & \mathrm{if}\,\ p\equiv -3 \pmod 8,\\
2p-1, & \mathrm{if}\,\ p\equiv 3 \pmod 8.
\end{array}
\right.
$$
 \end{thm}
Proof. With notations as before. That is, we use $S(X)$ the generating polynomial of $(s_u)_{u\ge 0}$
and let $\gamma\in GR^*(4^r,4)$ be of order $2p$ with  $S_0(\gamma)=\sum_{u \in D_0} \gamma^u \in GR^*(4^r,4)$.
Suppose that $C(X)\in \mathbb{Z}_4[X]$ is a connection polynomial of $(s_u)_{u\ge 0}$. We remark that $\min \deg(C(X))\le 2p$.

For  $p\equiv \pm 3 \pmod 8$, by \eqref{eq2} and Lemma \ref{l5}(II) we get
$$
C(\gamma^v)=0 ~~ \mathrm{for~ all}~~  v\in D_0\cup D_1\cup E_0\cup E_1.
$$
Now we consider the values of $C(\gamma^0)$ and $C(\gamma^p)$ \footnote{In fact, $C(\gamma^0)=C(1)$ and $C(\gamma^p)=C(-1)$.}.
Let $s(X)$ and $c(X)$ be the polynomials of degree $<2$ such that
$$
S(X)\equiv s(X) \pmod {X^2-1}, ~~ C(X)\equiv c(X) \pmod {X^2-1}.
$$

If $p\equiv -3 \pmod 8$, we have $S(-1)=S(1)=2$ by Lemma \ref{l5}(I).
It follows that $s(X)=2$ or $s(X)=2X$. So by (\ref{eq2}) again,
we see that
$$
c(X) \in \{0, 2, 2X, 2X+2\}
$$
and hence either
$C(-1)=C(1)=0$ or $C(-1)=C(1)=2$.

In terms of all values of $C(\gamma^v)$ for $v=0,1,\ldots,2p-1$
above, by Lemma \ref{root}(III) we have $\deg C(X)\geq 2p$.
Consequently, we get $\min \deg(C(X))= 2p$
  and hence $LC((s_u)_{u\ge 0})=2p$ for this case.

Similarly if $p\equiv 3 \pmod 8$, we have $S(1)=0$ and $S(-1)=2$ by Lemma \ref{l5}(I), and hence
$s(X)=1-X$. Then
we get
$$
c(X) \in \{0, 2, X+1, 2X+2\}
$$
and hence $C(-1)=C(1)=0$, or $C(-1)=C(1)=2$, or $C(-1)=0$ and
$C(1)=2$. Then by Lemma \ref{root}(III) we have $\deg C(X)\geq
2p-1$.

On  the other hand, since $s(X)=1-X$, we see that $S(X)$ is divisible by $X-1$ over $\mathbb{Z}_4$, from which we derive
$$
S(X)\cdot \frac{X^{2p}-1}{X-1} \equiv 0 \pmod {X^{2p}-1}.
 $$
Then $\frac{X^{2p}-1}{X-1}$  is a connection polynomial of $(s_u)_{u\ge 0}$.
So we get $\min \deg(C(X))= 2p-1$, i.e., $LC((s_u)_{u\ge 0})=2p-1$. \qed

\subsection{Linear complexity for the case $p\equiv \pm
1 \pmod 8$}
\label{subsec:LC-1mod8}

Due to Lemma \ref{l5}(III), it is more complicated to determine the connection polynomial $C(X)$ with the smallest degree
when
$p\equiv \pm 1 \pmod 8$. We need more detailed research on  the generating polynomial $S(X)$ of $(s_u)_{u\ge 0}$. For $j=0,1$, define
$$
\Gamma_j(X)=\prod_{v \in D_j}
(X- \gamma^v)  ~~\mathrm{and}~~ \Lambda_j(X)=\prod_{v \in E_j}
(X- \gamma^v),
$$
where $\gamma\in GR^*(4^r,4)$ is of order $2p$ with  $S_0(\gamma)=\sum_{u \in D_0} \gamma^u \in GR^*(4^r,4)$.
In particular, by Lemma \ref{relationD}(IV) we have for $j=0,1$,
$$
\Lambda_j(X)=\prod_{v \in D_j} (X-
\gamma^{v+p}).
$$

\begin{lem}\label{l6}
If  $p\equiv \pm 1 \pmod 8$,
then $\Gamma_j(X)$ and $\Lambda_j(X)$ are polynomials over $\mathbb{Z}_4$ for $j=0,1$.
\end{lem}
Proof. We only consider $\Gamma_0(X)$ here, for $\Gamma_1(X)$, $\Lambda_0(X)$ and $\Lambda_1(X)$ it can be done in a similar manner. It is sufficient to show
that the coefficients of $\Gamma_0(X)$
$$
a_m=(-1)^m\sum\limits_{\stackrel{i_1<i_2<\ldots<i_m}{i_1,i_2,\ldots,i_m \in D_0}}
\gamma^{i_1+i_2+\ldots+i_m}\in\mathbb{Z}_4
$$
for $1 \leq m \leq (p-1)/2$.

 Let $\gamma^{b}$ be a term of the last sum and $b \equiv i_1+i_2+\ldots+i_m \pmod
 {2p} $, $b \not \equiv 0\pmod
 {p}$. By Lemma \ref{relationD} for any $n: 0<n<(p-1)/2$ we have that $g^{2n}i_j \in D_0, j=0,\dots,m$.
  Hence, $X-\gamma^{g^{2n}i_j}, j=0,\dots,m$ are  the factors in the product $\prod_{v \in D_j}
(X- \gamma^v)$. So, $\gamma^{g^{2n}i_1}\dots
\gamma^{g^{2n}i_m}=\gamma^{g^{2n}b}$  is also a term of this sum
 for any
$n=0,...,(p-3)/2$, i.e.,
$$\gamma^{b}+\gamma^{g^2b}+\dots+\gamma^{g^{p-3}b}=S_0(\gamma^{b})$$
is a part of this sum. Therefore, there must exist the elements
$b_1, \dots,b_n$ such that
$$
 a_m=(-1)^m\sum\limits_{i=0}^n
S_0(\gamma^{b_i})+A,
$$
where  $A=|\{a| a\equiv (i_1+i_2+\ldots+i_m) \equiv 0 \pmod
{p}\mbox{ and }  i_1<i_2<\ldots<i_m; i_1,i_2,\ldots,i_m \in D_0\}|$.

By Lemma \ref{l4} and the proof of Lemma \ref{l1} we have that
$S_0(\gamma^{b_i}) \in\mathbb{Z}_4$. This  completes the proof of
Lemma \ref{l6}.   \qed

Since $\gamma^v$ is a root of $X^p+1$ for any $v\in \{p\}\cup D_0\cup D_1$, and
$\gamma^{v_1}-\gamma^{v_2}\in GR^*(4^r,4)$ for distinct $v_1,v_2\in \{p\}\cup D_0\cup D_1$,
it follows that
\begin{equation}\label{xp+1}
X^p+1=(X+1)\Gamma_0(X)\Gamma_1(X),
\end{equation}
from Lemma \ref{root} and the definitions on $\Gamma_0(X)$ and $\Gamma_1(X)$. Similarly, we have
\begin{equation}\label{xp-1}
X^p-1=(X-1)\Lambda_0(X)\Lambda_1(X).
\end{equation}
Now, let us explore the
expansion of
$$S(X)=2X^p+S_1(X)+2T_0(X)+3T_1(X).$$

\begin{lem}\label{l7}
We have the polynomial factoring in the ring $\mathbb{Z}_4[X]$
$$
\begin{array}{l}
S_1(X)+3T_1(X)=\\
\left\{
\begin{array}{cl}
(X^p-1)\Gamma_0(X)U_1(X),&\mathrm{if}\,\   p\equiv \pm 1 \pmod {16},\\
(X^p-1)\Gamma_0(X)U_2(X)+2(X^p+1),&\mathrm{if}\,\   p\equiv \pm 9 \pmod {16},
\end{array}
\right.
\end{array}
$$
and
$$
\begin{array}{l}
2X^p+2T_0(X)=\\
\left\{
\begin{array}{cl}
\Gamma_0(X)\Lambda_0(X)(X-1)V_1(X)+2(X^p+1),&\mathrm{if}\,\   p\equiv -1,-9 \pmod {16},\\
\Gamma_0(X)\Lambda_0(X)V_2(X)+2(X^p+1),&\mathrm{if}\,\   p\equiv 1,9 \pmod {16},
\end{array}
\right.
\end{array}
$$
where $U_i(X), V_i(X) \in \mathbb{Z}_4[X], i=1,2$.
\end{lem}
Proof.
Since $p\equiv \pm 1\pmod {8}$,  by Lemma \ref{relationD}(V) we obtain
\begin{eqnarray*}
S_1(X)+3T_1(X) &=& \sum_{u\in D_1}X^u+ 3\sum_{u\in E_1}X^u \\
& = & \sum_{\stackrel{u\in D_1}{u<p}}X^u+ \sum_{\stackrel{u\in D_1}{u>p}}X^u+3\sum_{\stackrel{u\in D_1}{u<p}}X^{u+p}+3\sum_{\stackrel{u\in D_1}{u>p}}X^{u-p}\\
&=&(X^p+3) \left(3\sum_{\stackrel{u\in D_1}{u<p}}X^{u}+\sum_{\stackrel{u\in D_1}{u>p}}X^{u-p}\right ).
\end{eqnarray*}
Write
$$M(X)=3\sum_{\stackrel{u\in D_1}{u<p}}X^{u}+\sum_{\stackrel{u\in D_1}{u>p}}X^{u-p}.$$
With the choice of $\gamma$ as before, if $v\in D_0$, we have
$$
M(\gamma^v)=3\sum_{\stackrel{u\in D_1}{u<p}}\gamma^{vu}+\sum_{\stackrel{u\in D_1}{u>p}}\gamma^{v(u-p)}=-\sum_{\stackrel{u\in D_1}{u<p}}\gamma^{vu}-\sum_{\stackrel{u\in D_1}{u>p}}\gamma^{vu}=
-S_1(\gamma)=-1+S_0(\gamma),
$$
where we use $\gamma^p=-1$ and (\ref{eq3}). So for $v\in D_0$, by Lemma \ref{l4} we get
 $$
 M(\gamma^v)=\left\{
\begin{array}{cl}
0,&\mathrm{if}\,\   p\equiv \pm 1 \pmod {16},\\
2,&\mathrm{if}\,\   p\equiv \pm 9 \pmod {16},
\end{array}
\right.
$$
from which, and by Lemma \ref{root}, we derive
$$
M(X)=\left\{
\begin{array}{cl}
\Gamma_0(X)U_1(X),&\mathrm{if}\,\   p\equiv \pm 1 \pmod {16},\\
2+\Gamma_0(X)U_2(X),&\mathrm{if}\,\   p\equiv \pm 9 \pmod {16},
\end{array}
\right.
$$
where $U_1(X), U_2(X) \in \mathbb{Z}_4[X]$.  We complete the proof of the first statement.

Now, we consider the polynomial $2X^p+2T_0(X)$. Since $2X^p+2T_0(X)=2(X^p+1)+2 +2 T_0(X)$, we only need to consider
$2 +2 T_0(X)$.

We first consider the roots of $2+2S_0(X)$. According to the proof of Lemma \ref{l1},
we see that $p+2\in D_0$ since $p\equiv \pm 1 \pmod {8}$. For any $v\in E_0$ with $v\equiv 2\overline{v} \pmod {2p}$, where $\overline{v}\in D_0$,
we obtain by (\ref{eq3}) and Lemma \ref{l4}
\begin{eqnarray*}
2+2S_0(\gamma^v)&=&2+2\sum_{u\in D_0} \gamma^{uv}=2+2\sum_{u\in D_0} \gamma^{2\overline{v}u}\\
&=&2+2S_0(\gamma^{2})=2+2S_0(-\gamma^{p+2})\\
&=&2-2S_0(\gamma)=0.
\end{eqnarray*}
So, by Lemma \ref{root} we
have $$2+2S_0(X)=\Lambda_0(X)G(X)$$
for some $G(X)\in \mathbb{Z}_4[X]$,  then we have
$$2+2S_0(X^2)=\Lambda_0(X^2)G(X^2).$$
Since $T_{0}(X)=S_0(X^2) \pmod{X^{2p}-1}$ and
\begin{eqnarray*}
 \Lambda_0(X^2)&=&\prod_{v
\in E_0} (X^2- \gamma^v)=\prod_{u \in D_0} (X^2-
\gamma^{2u})\\
&=&\prod_{u \in D_0} (X- \gamma^{u})(X+
\gamma^{u})\\
&=&\prod_{u \in D_0} (X- \gamma^{u})(X-
\gamma^{u+p})\\
&=&\prod_{u \in D_0} (X- \gamma^{u})\prod_{v \in D_0} (X-
\gamma^{v+p})\\
&=& \Gamma_0(X)\Lambda_0(X),
\end{eqnarray*}
it follows that
$$
2+2T_0(X)=2+2S_0(X^2)=\Gamma_0(X)\Lambda_0(X)G(X^2).
$$
On the other hand, from the fact that
$$
2+2T_0(1)=p+1=\left\{
\begin{array}{cl}
0,&\mathrm{if}\,\   p\equiv -1 \pmod {8},\\
2,&\mathrm{if}\,\   p\equiv 1 \pmod {8},
\end{array}
\right.
$$
and $\Gamma_0(1)\Lambda_0(1)\in GR^*(4^r,4)$, we write
$$
G(X^2)=(X-1)V_1(X)
$$
for $ p\equiv -1 \pmod {16}$ or $ p\equiv -9
\pmod {16}$. Otherwise, we write $V_2(X)=G(X^2)$.
Putting everything together, we complete the proof of the second statement. \qed

\begin{lem}\label{l8}
Let S(X) be the generating polynomial of $(s_u)_{u\ge
0}$ described in (\ref{SX}). We have in the ring $\mathbb{Z}_4[X]$
$$
S(X)=\left\{
\begin{array}{cl}
(X-1)\Gamma_0(X)\Gamma_1(X) W_1(X),&\mathrm{if}\,\   p\equiv -1 \pmod {16},\\
\Gamma_0(X)\Gamma_1(X)W_2(X),&\mathrm{if}\,\   p\equiv 1 \pmod {16},\\
(X-1)\Gamma_0(X)\Gamma_1(X) \Lambda_0(X)W_3(X),&\mathrm{if}\,\   p\equiv -9 \pmod {16},\\
\Gamma_0(X)\Gamma_1(X) \Lambda_0(X) W_4(X),&\mathrm{if}\,\   p\equiv
9 \pmod {16},
\end{array}
\right.
$$
where $W_i(\gamma^v)\neq 0, i=1,2$ for $v \in E_0\cup E_1$ and
$W_i(\gamma^v)\neq 0, i=3,4$ for $v \in  E_1.$
\end{lem}
Proof. Let $ p\equiv -1 \pmod {16}$. By (\ref{xp+1}) we have
$$2(X^p+1)=2(X+1)\Gamma_0(X)\Gamma_1(X)=2(X-1)\Gamma_0(X)\Gamma_1(X).$$
 Then according to
Lemma \ref{l7}, we write
$$
S(X)=(X-1)\Gamma_0(X)H(X),
$$
where $H(X)=U_1(X)(X^p-1)/(X-1)+\Lambda_0(X)V_1(X)+2\Gamma_1(X)$.
We check that
$$
H(\gamma^v)\left\{
\begin{array}{cl}
=0,&\mathrm{if}\,\   v\in D_1,\\
\neq 0,&\mathrm{if}\,\    v\in E_0,\\
\neq 0,&\mathrm{if}\,\    v\in E_1.
\end{array}
\right.
$$

 For $v\in D_1$, we have $S(\gamma^v)=0$ by Lemma \ref{l5}. Since $(\gamma^v-1)\Gamma_0(\gamma^v) \in  GR^*(4^r,4)$, we have $H(\gamma^v)=0$ by Lemma \ref{pre-root}.

For $v\in E_0$, we have $((\gamma^v)^p-1)/(\gamma^v-1)=0$ and $\Lambda_0(\gamma^v)=0$, so that $H(\gamma^v)=2\Gamma_1(\gamma^v)\neq 0$;

For $v\in E_1$, since $S(\gamma^v)= 2$ by Lemma \ref{l5}, we have $H(\gamma^v)\neq 0$.

\noindent
 So we have by Lemma \ref{pre-root}
$$
H(X)=\Gamma_1(X) W_1(X)
$$
for some $W_1(X) \in \mathbb{Z}_4[X]$ and $W_1(\gamma^v)\neq 0$ for $v \in E_0 \cup E_1$. Then we get
the factorization of $S(X)$ for $ p\equiv -1 \pmod {16}$.

Another assertions of this lemma can be
proved in a similar way. \qed

\begin{thm}
\label{t2}
Let $(s_u)_{u\ge 0}$ be the quaternary sequence over $\mathbb{Z}_4$ defined by \eqref{eq1}. Then the  linear complexity of $(s_u)_{u\ge 0}$
  satisfies
$$
LC((s_u)_{u\ge 0})
=\left\{
\begin{array}{cl}
p, & \mathrm{if}\,\ p\equiv -1 \pmod {16},\\
p+1, & \mathrm{if}\,\ p\equiv 1 \pmod {16}.
\end{array}
\right.
$$
\end{thm}
Proof. Let $p\equiv -1 \pmod {16}$. Since
$S(X)=(X-1)\Gamma_0(X)\Gamma_1(X) W_1(X)$ by Lemma \ref{l8}, together with Lemma \ref{l5}(I) we have
$W_1(\gamma^v)\neq 0$ for $v \in E_0 \cup E_1 \cup \{p\}$.
Then we see that
$$
S(X)(X+1)\Lambda_0(X)\Lambda_1(X) \equiv 0 \pmod {X^{2p}-1}.
 $$
That is, $(X+1)\Lambda_0(X)\Lambda_1(X)$ is a connection polynomial of degree $p$ of
$(s_u)_{u\ge 0}$. So the minimal degree of connection polynomials of
$(s_u)_{u\ge 0}$ is $\leq p$.

Let $C(X)\in \mathbb{Z}_4[X]$ be a connection polynomial of
$(s_u)_{u\ge 0}$. Due to
$$
\gcd( (X-1)\Gamma_0(X)\Gamma_1(X), (X+1)\Lambda_0(X)\Lambda_1(X) )=1,
$$
we have
$$
W_1(X)C(X)\equiv 0 \pmod
{(X+1)\Lambda_0(X)\Lambda_1(X)}
$$
by \eqref{eq2}, \eqref{xp+1}, \eqref{xp-1}  and Lemma \ref{l8}.
So we deduce
$$W_1(\gamma^v)C(\gamma^v)=0 ~~ \mathrm{for} ~~ v \in E_0 \cup E_1 \cup \{p\}.$$

Since $W_1(\gamma^v)\neq
0$ for $v \in E_0 \cup E_1 \cup \{p\}$, if $W_1(\gamma^v)
\in GR^*(4^r,4)$ then we get $C(\gamma^v)= 0$, and if $W_1(\gamma^v)=2\eta, \eta \in GR^*(4^r,4)$,
then we have either $C(\gamma^v)= 0$ or $C(\gamma^v)= 2$. I.e., $2C(\gamma^v)= 0$  for $v \in E_0 \cup E_1 \cup
\{p\}.$

By the definition of $C(X)=1+c_1X+\ldots$, i.e., $2C(x)$ is non-constant, then by
Lemma \ref{pre-root} we have that $2C(x)$ is
 divisible by $(X+1)\Lambda_0(X)\Lambda_1(X)$, i.e.,
$\deg C(X)\geq p$ and hence $LC((s_u)_{u\ge 0})=p$ for this case. We prove the first statement.

Let $p\equiv 1 \pmod {16}$. From
that
$$
S(X)(X^2-1)\Lambda_0(X)\Lambda_1(X) \equiv 0 \pmod {(X^{2p}-1)},
 $$
we see that $(X^2-1)\Lambda_0(X)\Lambda_1(X)$ is a connection polynomial of
$(s_u)_{u\ge 0}$  of degree $p+1$.

For any connection polynomial $C(X)$ of $(s_u)_{u\ge 0}$, a similar way presented above gives
$$W_2(X)C(X)\equiv 0 \pmod
{(X^2-1)\Lambda_0(X)\Lambda_1(X)}. $$
 As in the proof of Theorem \ref{t1}, denote by $s(X)$ and $c(X)$  the polynomials of degree $<2$ such
that
$$
S(X)\equiv s(X) \pmod {X^2-1}, ~~ C(X)\equiv c(X) \pmod {X^2-1}.
$$
As earlier, we  can obtain  that $ c(X) \in \{0, 2, 2X, 2X+2\} $,
hence $2c(X)=0$ and $2C(X)=(X^2-1)M(X)$ for some $M(X)\in \mathbb{Z}_4[X]$. Since by Lemma \ref{l8}
$W_2(\gamma^v)\neq 0$ for $v \in E_0 \cup E_1 $, it follows that
$2C(\gamma^v)= 0$ and $M(\gamma^v)=0$ for $v \in E_0 \cup E_1.$ Therefore,
$M(X)$ is
 divisible by $\Lambda_0(X)\Lambda_1(X)$ by Lemma \ref{pre-root}, i.e.,
$\deg C(X)\geq p+1$ and hence $LC((s_u)_{u\ge 0})=p+1$ for this
case.
  \qed

\begin{thm}
\label{t3}
Let $(s_u)_{u\ge 0}$ be the quaternary sequence defined by \eqref{eq1}. Then the  linear complexity of $(s_u)_{u\ge 0}$
  satisfies
$$
LC((s_u)_{u\ge 0})
=\left\{
\begin{array}{cl}
(p+1)/2, & \mathrm{if}\,\ p\equiv -9 \pmod {16},\\
(p+3)/2, & \mathrm{if}\,\ p\equiv 9 \pmod {16}.
\end{array}
\right.
$$
\end{thm}
Proof. The proof can follow that of Theorem \ref{t2} in a similar way. Here we give a sketch.

Let $p\equiv -9 \pmod {16}$.  On the one hand,
$(X+1)\Lambda_1(X)$ is a connection polynomial of
$(s_u)_{u\ge 0}$ of degree $(p+1)/2$ by \eqref{eq2}.

On the other hand, for any connection polynomial $C(X)$ of $(s_u)_{u\ge 0}$, we have
$$W_3(X)C(X)\equiv 0 \pmod {(X+1)\Lambda_1(X)}.$$
Now since $W_3(\gamma^v)\neq 0$ for $v \in E_1 \cup \{p\}$, it follows
that $2C(\gamma^v)= 0$ for $v \in \cup E_1 \cup \{p\}.$ Therefore,
by Lemma \ref{root} again $2C(x)$ is
 divisible by $(X+1)\Lambda_1(X)$, i.e.,
$\deg C(X)\geq (p+1)/2 $ and hence $LC((s_u)_{u\ge 0})=(p+1)/2$ for
this case.

The case of $p\equiv 9 \pmod {16}$ follows the way of $p\equiv 1 \pmod {16}$ in Theorem \ref{t2} and we omit it. \qed

\section{Final remarks and conclusions}
\label{sec:3}

We determined the exact values of the linear complexity of $2p$-periodic quaternary
sequences over $\mathbb{Z}_4$ defined from the generalized cyclotomic
classes modulo $2p$ by considering the factorization of the generating polynomial $S(X)$ in $\mathbb{Z}_4[X]$. It is more complicated to study this problem than that in finite fields. Besides the autocorrelation considered in \cite{KHS}, this is another cryptographic feature of the quaternary cyclotomic sequences of period $2p$.

A direct computing of the linear complexity has been done for $3
\leq p \leq 1000$ by the Berlekamp-Massey algorithm adapted by Reeds
and Sloane  in \cite{RS} for the residue class ring  to confirm our
theorems. Below we list some experimental data.

1. $p=3$, $(s_u)_{u\ge 0}=(0,0,2,2,3,1)$, then
$C(X)=1+X+X^2+X^3+X^4+X^5$ and $LC((s_u)_{u\ge 0})=5(=2p-1).$

2. $p=5$, $(s_u)_{u\ge 0}=(0,0,2,1,3,2,3,1,2,0)$, then
$C(X)=1+3X^{10}$ and $LC((s_u)_{u\ge 0})=10(=2p).$

3. $p=7$, $(s_u)_{u\ge 0}=(0,0,2,1,2,1,3,2,2,0,3,0,3,1)$, then
$C(X)=1+X^2+X^3+3X^4$ and $LC((s_u)_{u\ge 0})=4(=(p+1)/2).$

4. $p=17$, $C(X)=1+X+3X^{17}+3X^{18}$, $LC((s_u)_{u\ge 0})=18(=p+1).$

5. $p=31$, $C(X)=1+3X^{31}$, $LC((s_u)_{u\ge 0})=31(=p).$

6. $p=41$,
$C(X)=1+2X^2+3X^3+2X^5+2X^6+3X^7+3X^8+3X^9+X^{10}+2X^{11}+3X^{12}+X^{13}+X^{14}+X^{15}+2X^{16}+2X^{17}+X^{19}+2X^{20}+3X^{22}$,
$LC((s_u)_{u\ge 0})=22(=(p+3)/2).$

We hope that the procedures in this paper used to
derive the linear complexity can be extended to quaternary cyclotomic sequences with
larger period (for example, $2p^n$).

We finally remark that it is interesting to consider the $k$-error linear complexity of the sequences in this work.

\section*{Acknowledgements}
Parts of this work were written during a very pleasant visit of Z. Chen to the Hong Kong University of Science and Technology. He wishes to thank Prof. Cunsheng Ding for the hospitality and financial support.

Z.X.C. was partially supported by the National Natural Science
Foundation of China under grant No. 61373140.

 V.A.E. was supported by the Ministry of Education and Science
of Russia as a part of state-sponsored project No. 1.949.2014/K.

\section*{Appendix}

\emph{Proof of Lemma 1}.

It is a well-known fact that if $\xi \in
GR(4^r,4)$ is a root of the polynomial $P(X)\in \mathbb{Z}_4[X]$
then $P(X)=(X-\xi)Q_1(X)$ for some polynomial  $Q_1(X)\in
GR(4^r,4)[X]$. Let $\eta$ be another root of   $P(X)$ and $\xi-\eta
\in GR^*(4^r,4)$ then $0=(\xi-\eta) Q_1(\eta)$, i.e.,$Q_1(\eta)=0$.
So, $Q_1(X)=(X-\eta)Q_2(X)$ for some polynomial $Q_2(X)\in
GR(4^r,4)[X]$ and $P(X)=(X-\xi)(X-\eta)Q_2(X)$.\\

\emph{Proof of Lemma 2.}

 (I). By the choice of $\gamma$  we have an
expansion $(X^p-1)/(X- 1)=\prod_{j=1}^{p-1}(X- \gamma^{2j})$, hence
$p =\prod_{j=1}^{p-1}(1-\gamma^j)(1+\gamma^j)$ . So,
$\gamma^j-\gamma^n \in GR^*(4^r,4)$  when $j, n=0,...,2p-1$ and $j
\not \equiv n~(\bmod~p)$. Therefore, if $P(\gamma^j)=0$ for all $j
\in D_i$ or for all $j \in E_i, i=0,1$, then by Lemma \ref{pre-root}
$ P(X)$ is divisible by $ \prod_{j \in D_i} (X-\gamma^j)$  or  $
\prod_{j \in E_i} (X-\gamma^j)$.  The first assertion of Lemma
\ref{root} is proved.

(II). This assertion follows from (I).

 (III). We consider two cases.

 Let $P(1)=0$ or $P(-1)=0$. Suppose $P(-1)=0$, in this case by (II) we have that  $$P(X)=\left
(X^p+1 \right)P_3(X)$$ and $2P_3(X)\neq 0$ since $P(0)=1$. From the
equality $P(X)=\left (X^p+1 \right)P_3(X)$ for $X=\gamma^v,$ $v \in
E_0 \cup E_1 $  we deduce $2P_3(\gamma^v)=0$, therefore $2P_3(X)$ is
divisible by $(X^p-1)/(X-1)$ and $\deg P(X) \geq 2p-1$. Furthermore,
if $P(1)=0$ then $2P_3(X)$ is divisible by $(X^p-1)$ and $\deg P(X)
\geq 2p$.

 Let $P(1)\neq 0$ and $P(-1)\neq 0$. Then, by condition  $P(1)=2$ and
 $P(-1)=2$ and $P(\gamma^j)=0$ for all $j: j\in D_0 \cup D_1$. By (I) we
have that   $P(X)=Q(X)\left (X^p+1 \right)/(X+1)$, $Q(X) \in
\mathbb{Z}_4 [X]$ and $2Q(X)\neq 0$. Since $P(-1)=2$ it follows that
$Q(-1)=2$ and $Q(X)=(X+1)F(X)+2$, $F(X) \in \mathbb{Z}_4 [X]$ or
 $$P(X)=\left (X^p+1
\right)F(X)+2\left (X^p+1 \right)/(X+1).$$ From the last equality
and conditions of this lemma we obtain that $2F(\gamma^v)=0$  for $v
\in E_0 \cup E_1 \cup \{0\}$, therefore $2F(x)$ is divisible by
$x^p-1$ and $\deg P(X)\geq 2p$.

\emph{Remark.} The polynomial $P(X)$ is not obliged to be  divisible
by $X^{2p}-1$ when $P(\gamma^j)=0$ for $j=0,1,\dots,2p-1$. For
example, $P(X)= X^{2p}-1 +2(X^p+1)$.\\


\begin{thebibliography}{12}


\bibitem{CL}
Z. L. Chang, D. D. Li.
On the linear complexity of quaternary cyclotomic sequences with the period $2pq$.
 IEICE Trans. Fundamentals of Electronics, Communications
and Computer Sciences, 2014, E97-A(2): 679-684.

\bibitem{CDR} T. W. Cusick, C. Ding, A. Renvall. Stream Ciphers and Number
 Theory. Elsevier. Amsterdam (1998).

\bibitem{DHM} C. Ding,  T.Helleseth, H. Martinsen. New families of binary sequences with optimal three-level autocorrelation.
IEEE Transactions on Information Theory, 2001, 47(1): 428-433.


\bibitem{DC} X. Du, Z. Chen. Linear complexity of quaternary sequences
generated using generalized cyclotomic classes modulo $2p$.
  IEICE Trans. Fundamentals of Electronics, Communications
and Computer Sciences, 2011, E94-A(5): 1214-1217.


\bibitem{E} V.A. Edemskiy. On the linear complexity of binary sequences on the basis of biquadratic and sextic residue
classes. Discret. Math. Appl., 2010, 20(1): 75-84 (2010) (Diskretn.
Mat., 2010, 22(1): 74-82).


\bibitem{EI1}
V. Edemskiy, A. Ivanov.  Autocorrelation and linear complexity of
quaternary sequences of period $2p$ based on cyclotomic classes of
order four. Proceedings of the 2013 IEEE International Symposium on
Information Theory, ISIT 2013. July 7-12,  Istanbul, Turkey,
3120-3124.

\bibitem{EI2}
V. Edemskiy, A. Ivanov. Linear complexity of quaternary sequences of
length $pq$ with low autocorrelation. Journal of Computational and
Applied Mathematics, 2014, 259 : 555-560.

\bibitem{EI3}
V. Edemskiy, A. Ivanov. The linear complexity of balanced quaternary
sequences with optimal autocorrelation value. Cryptography and
Communications, 2015, 7(4): 485-496.

 \bibitem{GG} S.W. Golomb, G. Gong. Signal Design for Good Correlation:
ForWireless Communications, Cryptography and Radar Applications.
Cambridge University Press, Cambridge (2005).

 \bibitem{H} M. Hall.
Combinatorial Theory. Wiley, New York (1975).

\bibitem{IR} K. Ireland, M.Rosen. A Classical Introduction to Modern Number
Theory. Springer (1982).

\bibitem{KYZ}
P. Ke,  Z. Yang,  J. Zhang. On the autocorrelation and linear
complexity of some 2p periodic quaternary cyclotomic sequences over
$\mathbb{F}_4$.  IEICE Trans. Fundamentals of Electronics, Communications and
Computer Sciences, 2011,  EA94-A(11): 2472-2477.

\bibitem{KZ}
P. Ke, S. Zhang. New classes of quaternary cyclotomic
sequence of length $2p^m$ with high linear complexity. Information Processing
Letters, 2012, 112: 646-650.


\bibitem{KHS}  Y.-J. Kim, Y.-P. Hong, H.-Y. Song. Autocorrelation of some quaternary cyclotomic sequences of length
$2p$. IEICE Trans. Fundamentals of Electronics, Communications and
Computer Sciences, 2008, E91-A(12): 3679-3684.



\bibitem{M} B. R. McDonald. Finite Rings With Identity. New York. Marcel
Dekker (1974).

\bibitem{RS} J. A. Reeds, N. J. A. Sloane. Shift-register synthesis modulo $m$. SIAM J. Comput., 1968,  14:
505-513.




\bibitem{US96}
P. Udaya, M. U. Siddiqi. Optimal biphase sequences with large linear complexity derived from sequences over $\mathbb{Z}_{4}$.
IEEE Transactions on Information Theory, 1996, 42 : 206-216.

\bibitem{US98}
P. Udaya, M. U. Siddiqi. Optimal and suboptimal quadriphase sequences derived from maximal length sequences
over $\mathbb{Z}_{4}$. Applicable Algebra in Engineering, Communication and Computing, 1998, 9 : 161-191.


\bibitem{US00}
P. Udaya, M. U. Siddiqi.
    Generalized GMW Quadriphase Sequences
Satisfying the Welch Bound with Equality. Applicable Algebra in Engineering, Communication and Computing,
 2000, 10 : 203-225.




\bibitem{W1}
 Z.X. Wan. Finite Fields and Galois Rings. Singapore. World
Scientific Publisher (2003).




\bibitem{ZZ15}
J. Zhang, C. A. Zhao. The linear complexity of a class of binary sequences with period $2p$. Applicable Algebra in Engineering,
Communication and Computing, 2015, 26(5): 475-491.


\end{thebibliography}
\end{document}